\documentclass{amsart}

\usepackage{amssymb}
\usepackage{eepic}
\usepackage{color}
\usepackage[linktocpage=true]{hyperref}

\allowdisplaybreaks

\numberwithin{equation}{section}

\newtheorem{theorem}{Theorem}[section]
\newtheorem{lemma}[theorem]{Lemma}

\newtheorem{remark}[theorem]{Remark}

\newtheorem{TheoA}{Theorem A}

\newtheorem{TheoB}{Theorem B}

\newcommand{\summ}{\sum\nolimits}

\def\1{\mathbf{1}}
\def\H{\mathcal{H}}
\def\E{\mathcal{E}}

\def\M{\mathcal{M}}

\newcommand{\dem}{\noindent {\bf Proof. }}

\newcommand{\fin}{\hspace*{\fill} $\square$ \vskip0.2cm}

\addtocounter{tocdepth}{0}

\begin{document}

\null

\vskip-40pt

\null

\title[Asymmetric Doob inequalities]{Asymmetric Doob inequalities \\ in continuous time}

\author[Hong, Junge, Parcet]
{Guixiang Hong, Marius Junge, Javier Parcet}

\maketitle

\null

\vskip-45pt

\null

\maketitle

\begin{abstract}
The present paper is devoted to the second part of our project on asymmetric maximal inequalities, where we consider martingales in continuous time. Let $(\M,\tau)$ be a noncommutative probability space equipped with a continuous filtration of von Neumann subalgebras $(\M_t)_{0\leq t\leq1}$ whose union is weak-$*$ dense in $\mathcal{M}$. Let $\E_t$ denote the corresponding family of conditional expectations. As for discrete filtrations, we shall prove that for $1 < p < 2$ and $x \in L_p(\M,\tau)$ one can find $a, b \in L_p(\M,\tau)$ and contractions $u_t, v_t \in \M$ such that $$\E_t(x) = a u_t + v_t b \quad \mbox{and} \quad \max \big\{ \|a\|_p, \|b\|_p \big\} \le c_p \|x\|_p.$$ Moreover, $a u_t$ and $v_t b$ converge in the row/column Hardy spaces $\H_p^r(\M)$ and $\H_p^c(\M)$ respectively. We also confirm in the continuous setting the validity of related asymmetric maximal inequalities which we recently found for discrete filtrations, including $p=1$. As for other results in noncommutative martingale theory, the passage from discrete to continuous index is quite technical and requires genuinely new methods. Our approach towards asymmetric maximal inequalities is based on certain construction of conditional expectations for a sequence of projective systems of $L_p$-modules. The convergence in $\H_p^r(\M)$ and $\H_p^c(\M)$ also imposes new algebraic atomic decompositions.
\end{abstract}

\section*{{\bf Introduction}}

Noncommutative martingales associated to continuous filtrations arise naturally in the dilation theory of Markov semigroups over von Neumann algebras \cite{JES}. In this direction, a theory of Hardy spaces for continuous filtrations was formulated in \cite{JuPe14,Per}. This includes (conditional) Hardy and BMO spaces for noncommutative continuous martingales together with Fefferman-Stein duality, Burkholder-Gundy inequalities, Davis and Burkholder-Rosenthal inequalities, etc. The theory admits applications in quantum stochastic calculus and especially in harmonic analysis through the connection with the interpolation theory of BMO spaces associated with semigroups of operators \cite{JM2}, which plays a crucial role in understanding the behavior in $L_p$ of noncommutative Riesz transforms \cite{JM1} and other smooth Fourier multipliers in group von Neumann algebras \cite{JMP}. 

In our recent paper \cite{HJP}, we established new asymmetric Doob inequalities based on an algebraic form of the atomic decomposition which yield a very satisfactory and complete picture of the noncommutative maximal inequalities for martingales in discrete time. This included the right reformulation of Davis' martingale theorem \cite{Da} for $p=1$, which escaped previous attempts for quite some time. It is therefore very natural to wonder whether similar estimates hold in continuous time. Beyond its unquestionable interest in the theory of Hardy spaces, this could lay foundations for future applications in maximal ergodic theory and harmonic analysis. As we shall see in this paper, asymmetric maximal inequalities in continuous time do not follow from the discrete case via somehow standard limiting processes, but implementing new techniques involving $L_p$-modules and introducing new families of Hardy spaces to generalize our original arguments. 

Given $(\M,\tau)$ a noncommutative probability space, let $\E_t$ denote the conditional expectations associated to a weak-$*$ dense filtration $(\M_t)_{t\in[0,1]}$ of von Neumann subalgebras. Let $H^c_p(\sigma)$ be the column Hardy space with respect to the discrete filtration $(\M_t)_{t\in\sigma}$ for a finite partition $\sigma=\{0=t_0<t_1<\dotsm t_n=1\}$ of the interval $[0,1]$. Fix an ultrafilter $\mathcal{U}$ over the set of all such finite partitions $\Sigma$ such that for each finite partition $\sigma$ of $[0,1]$ the set $U_\sigma=\{\sigma'\in\Sigma:\sigma\subset\sigma'\}$ belongs to $\mathcal{U}$. Given $x\in\mathcal{M}$ and $1\leq p<\infty$, set
$$\|x\|_{\mathcal{H}^c_p}=\lim_{\sigma,\mathcal{U}}\|x\|_{H^c_p(\sigma)}.$$ The Hardy space $\mathcal{H}^c_p(\M)$ associated to the continuous filtration $(\M_t)_{t\in[0,1]}$ is the completion of $\M$ with respect to the norm $\|\cdot\|_{\mathcal{H}^c_p}$ which |based on the monotonicity properties of the $H^c_p(\sigma)$-norm| was proved to be independent of the ultrafilter $\mathcal{U}$ in \cite{JuPe14}. For this reason, we will work with one fixed ultrafilter $\mathcal{U}$ over $\Sigma$ in the present paper. Row Hardy spaces are defined analogously. 

The asymmetric maximal quasi-norms for continuous indices are trivially defined as in the discrete case, we quickly introduce them here and refer to \cite{HJP} for further information. Given $1\leq p<\infty$ and $0\leq\theta\leq1$, we take $$\big\| (x_t)_{t\in[0,1]} \big\|_{L_p(\ell_\infty^\theta)} = \inf \Big\{ \|a\|_{\frac{p}{1-\theta}} \Big( \sup_{t\in[0,1]} \|w_t\|_\infty \Big) \|b\|_{\frac{p}{\theta}} \, \big| \ x_t= a w_t b \mbox{ for } t\in[0,1]\Big\}.$$
Additionally, the weak column maximal quasi-norm is defined as follows $$\big\| (x_t)_{t\in[0,1]} \big\|_{\Lambda_{p, \infty}(\ell_\infty^c)} \hskip1pt = \hskip1pt \sup_{\lambda > 0} \inf_{q \in \M_\pi} \Big\{ \lambda \tau ( \1 - q )^{\frac{1}{p}} \, \big| \ \| x_t q\|_\infty \le \lambda \mbox{ for all } t\in[0,1] \Big\}.$$ $\M_\pi$ denotes the projection lattice in $\M$. Take adjoints to define row quasi-norms. 

\begin{TheoA}
Let $(\M,\tau)$ be a noncommutative probability space and let $\E_t$ denote the conditional expectations associated to a weak-$*$ dense filtration $(\M_t)_{t\in[0,1]}$ of von Neumann subalgebras. Then, the following inequalities hold$\hskip1pt :$
\begin{itemize}
\item[i)] Given $1 \le p \le 2$ and $x \in \H_p^c(\M)$ $$\big\| (\E_t(x))_{t\in[0,1]} \big\|_{\Lambda_{p,\infty}(\ell_\infty^c)} \le c_p \|x\|_{\H_p^c}.$$
The row analog $(\E_t)_{t\in[0,1]}: \H_p^r(\M) \to \Lambda_{p,\infty}(\M;\ell_\infty^r)$ is also bounded.

\vskip3pt

\item[ii)] Given $1 \le p \le 2$ and $x \in \H_p^c(\M)$ $$\big\| (\E_t(x))_{t\in[0,1]} \big\|_{L_p(\ell_\infty^\theta)} \le c_{p,\theta} \|x\|_{\H_p^c}$$ provided $1 - p/2 < \theta < 1$. The same holds for $x \in \H_p^r(\M)$ and $0 < \theta < p/2$.
\end{itemize}
\end{TheoA}

Theorem Ai will be deduced directly from the result in the discrete case since ultralimit of projections can be easily modified to be a projection satisfying the desired properties.
A crucial difficulty in the proof of Theorem Aii is that we may not use the identity $$\| (\E_t(x))_{t\in[0,1]} \big\|_{L_p(\M;\ell_\infty^\theta)}=\sup_{\sigma\in\Sigma} \big\| (\E_t(x))_{t\in\sigma} \big\|_{L_p(\M;\ell_\infty^\theta)}$$ when $p/2<\theta<1$, since $\|\cdot\|_{L_p(\ell_\infty^\theta)}$ is not a norm. Thus, the discrete time results in \cite{HJP} can not be used as a black box. Instead, we give a direct argument. A key new ingredient is the construction of a sequence of increasing von Neumann algebras from a sequence of projective systems of $L_p$ modules, see Lemma \ref{lem:construction of conditional expectations}.

Stronger asymmetric Doob maximal estimates follow by stretching our approach to produce finer algebraic Davis type decompositions. More precisely, motivated by the proof of Theorem A, we will introduce new Hardy spaces $$\underbrace{\widehat{h}_{pw}^r(\M) + \widehat{h}_{pw}^{1_r}(\M)}_{\widehat{\H}_{pw}^r(\M)} \quad \mbox{and} \quad \underbrace{\widehat{h}_{pw}^c(\M) + \widehat{h}_{pw}^{1_c}(\M)}_{\widehat{\H}_{pw}^c} \qquad \mbox{for $w \ge 2$}.$$ 

\begin{TheoB}
Let $(\M,\tau)$ be a noncommutative probability space and let $\E_t$ denote the conditional expectations associated to a weak-$*$ dense filtration $(\M_t)_{t\in[0,1]}$ of von Neumann subalgebras. Then, the following results hold$\hskip1pt :$
\begin{itemize}
\item[i)] Given $1 < p < 2$ with $1/p = 1/w + 1/s$, we find $$L_p(\M) \simeq \widehat{\H}_{pw}^r(\M) + \widehat{\H}_{pw}^c(\M) \quad \mbox{provided} \quad w>2, s \ge 2.$$

\vskip3pt

\item[ii)] Given $1 < p < 2$, the inequalities below hold for any $w > 2$
\begin{eqnarray*}
\big\| (\E_t(x))_{t\in[0,1]} \big\|_{L_{p}(\ell_\infty^r)} \!\!\!\! & \le & \!\!\!\! c_{pw} \|x\|_{\widehat{\H}_{pw}^r}, \\
\big\| (\E_t(x))_{t\in[0,1]} \big\|_{L_{p}(\ell_\infty^c)} \!\!\!\! & \le & \!\!\!\! c_{p,w} \|x\|_{\widehat{\H}_{pw}^c}.
\end{eqnarray*}
\noindent In particular, given $x \in L_p(\M)$ may write $x = x_r + x_c$ with $$\hskip30pt \max \Big\{ \big\| (\E_t(x_r))_{t\in[0,1]} \big\|_{L_{p}(\ell_\infty^r)}, \big\| (\E_t(x_c))_{t\in[0,1]} \big\|_{L_{p}(\ell_\infty^c)} \Big\} \, \le \, c_p \|x\|_{p}.$$

\vskip3pt

\item[iii)] Moreover, we have $x_r \in \widehat{\H}_{pw}^r(\M) \subset \H_p^r(\M)$ and $x_c \in \widehat{\H}_{pw}^c(\M) \subset \H_p^c(\M)$.
\end{itemize}
\end{TheoB}

Theorem Biii establishes an important fact. Namely, that the decomposition elements can be taken in the usual Hardy spaces defined via square functions. To that end, we prove an algebraic Davis decomposition for these newly defined Hardy spaces in Theorem \ref{pro:hat hcpw=hcpw}. We think this is of independent interest. 

\section{{\bf Proof of Theorem A}}

In this section we establish our first asymmetric inequalities in continuous time already stated in Theorem A. Unlike the symmetric maximal inequalities, which follow trivially from the discrete results, the asymmetric ones in Theorem Aii require ultraproduct methods and $L_p$ module theory to generalize the algebraic atomic decompositions of $L_p$ and Hardy spaces.

\subsection{Proof of Theorem Ai} This is relatively easy since we can reduce it to the analogous result in the discrete case \cite{HJP}. Indeed, by density we may assume that $x \in \M$. Now, given $\lambda>0$ and $\sigma\in \Sigma$, we know from \cite[Theorem Aii]{HJP} that there exists a projection $q_\sigma\in \M$ such that 
$$\big\| \mathcal{E}_t(x) q_\sigma \big\|_\infty \le \lambda \quad \mathrm{and} \quad \lambda \tau\big( \1 - q_\sigma \big)^{\frac1p} \le C_p \|x\|_{{H}^c_p(\sigma)}$$ for any $t \in\sigma$. Define $u=w^*-L_{\infty}-\lim_{\sigma,\mathcal{U}}q_\sigma$. Recall that $u$ is not necessarily a projection. However, recalling that $x \in \M$, it is straightforward to show that the exact same inequalities above apply for $u$ instead of $q_\sigma$ and $\H_p^c(\M)$ instead of $H_p^c(\sigma)$, details are left to the reader. Then, the projection $q = \chi_{[\frac12,1]}(u)$ clearly satisfies $$q \le 2quq \le 4 u^2 \quad \mbox{and} \quad \1-q \le 2(\1-u).$$ This implies for any $t \in[0,1]$ that \\ [3pt] \null \hfill \hskip5pt 
$\big\|\mathcal{E}_t(x)q \big\|_\infty \le 2 \big\| \mathcal{E}_t(x) u^2 \mathcal{E}_t(x)^* \big\|_\infty^\frac12 \le 2 \lambda \quad \mbox{and} \quad \lambda \tau \big( \1 - q \big)^{\frac1p} \le 2^{\frac1p} C_p \|x\|_{\mathcal{H}^c_p}.$ \hfill $\square$

\subsection{Ultraproducts and projective systems of $L_p$ modules}

Before proving Theorem Aii let us recall some preliminaries, we refer to \cite{JuPe14} and references therein for more information. Let $\mathcal{U}$ be the ultrafilter that we have fixed over $\Sigma$ and consider a family of Banach spaces $(X^\sigma)_{\sigma\in\Sigma}$. Let $\ell_\infty(\{X^\sigma : \sigma \in\Sigma\})$ be the space of bounded families $(x^\sigma)_{\sigma\in\Sigma}\in\prod_{\sigma}X^\sigma$ and define the ultraproduct $\prod_{\mathcal U} X^\sigma$ as the quotient space $$\prod_\mathcal{U} X^\sigma = \ell_\infty \big( \{X^\sigma : \sigma\in\Sigma\} \big) \big/ \mathcal{N}^{\mathcal U},$$ where $\mathcal N^\mathcal U$ denotes the (closed) subspace of $\mathcal U$-vanishing families $$\mathcal N^\mathcal U = \Big\{ (x^\sigma)_{\sigma\in\Sigma} \in\ell_\infty(\{X^\sigma : \sigma \in\Sigma \}) \, \big| \, \lim_{\sigma,\mathcal U}\|x^\sigma\|_{X^\sigma} = 0 \Big\}.$$ If $(x^\sigma)^\bullet$ is the element of $\prod_\mathcal U X^\sigma$ represented by $(x^\sigma)_{\sigma\in\Sigma}$, its quotient norm is 
$$\|(x^\sigma)^\bullet\| = \lim_{\sigma, \mathcal U}\|x^\sigma\|_{X^\sigma}.$$ When $X^\sigma = X$ for all $\sigma$, the ultrapower $\prod_\mathcal U X$ is the quotient space $\ell_\infty(\Sigma;X)/\mathcal N^\mathcal U$. If $(X^\sigma)_{\sigma\in\Sigma}$ and $(Y^\sigma)_{\sigma\in\Sigma}$ are two families of Banach spaces and $T^\sigma : X^\sigma \rightarrow Y^\sigma$ are linear operators uniformly bounded in $\sigma \in\Sigma $, the ultraproduct map $(T^\sigma)^\bullet$ is defined canonically as follows $$(T^\sigma)^\bullet:\prod_\mathcal U X^\sigma\ni (x^\sigma)^\bullet \mapsto (T^\sigma x^\sigma)^\bullet\in \prod_\mathcal U Y^\sigma.$$ We refer to \cite{Hei80,Sim82} for basic facts about ultraproducts of Banach spaces. 

It is well-known that the class of von Neumann algebras is not closed under ultrapowers, but according to the work of Groh \cite{Gro84}, the class of the preduals of von Neumann algebras is. Let $\mathcal M$ be a von Neumann algebra. Then $\prod_\mathcal U \mathcal M_*$ is the predual of a von Neumann algebra denoted by
$$\widetilde{\mathcal M}_\mathcal U = \Big( \prod_\mathcal U \mathcal M_* \Big)^*.$$ $\prod_\mathcal U\mathcal M$ may be identified with a weak*-dense subalgebra and $\widetilde{\mathcal M}_\mathcal U$ becomes the von Neumann algebra generated by $\prod_\mathcal U \mathcal M$ in $\mathcal B(\prod_\mathcal U H)$ for any $*$-representation of $\mathcal M$ in $\mathcal{B}(H)$. This was proved by Raynaud \cite{Ray02}, who also constructed an isometric isomorphism
$$L_p (\widetilde{\mathcal M}_\mathcal U) \simeq \prod_\mathcal U L_p ({\mathcal M}) \quad \mbox{for all $p > 0$}.$$ If $\M$ is a finite, the usual von Neumann algebra ultrapower is $$\mathcal M_\mathcal U = \ell_\infty(\Sigma;\mathcal M)/\Sigma^\mathcal U \ \ \mbox{where} \ \ \Sigma^\mathcal U = \Big\{(x^\sigma)_{\sigma\in\Sigma}\in\ell_\infty(\Sigma; \mathcal M) \, \big| \, \lim_{\sigma,\mathcal U} \tau((x^\sigma)^* x^\sigma) = 0 \Big\}.$$ According to Sakai \cite{Sak71}, $\mathcal M_\mathcal U$ is a finite von Neumann algebra when equipped with the ultrapower map of the trace $\tau_\mathcal U ((x^\sigma)^\bullet) = \lim_{\sigma,\mathcal U} \tau(x^\sigma). $
This is compatible with $\Sigma^\mathcal U$ and defines a normal faithful normalized trace on $\mathcal M_\mathcal U$. We may identify $\mathcal M_\mathcal U$ as a dense subspace of $L_1(\mathcal M_\mathcal U)$ via the map $\mathcal M_\mathcal U \ni x \mapsto \tau_\mathcal U(x \hskip1pt \cdot)\in L_1(\mathcal M_\mathcal U)$. Then we have $\|x\|_1 = \lim_{\sigma,\mathcal U}\|x^\sigma\|_1$ for $x = (x^\sigma)^\bullet\in\mathcal M_\mathcal U$, regardless of which representing family $(x^\sigma)$ of $x$ we use. This implies that $L_1(\M_\mathcal{U})$ can be regarded as an isometric subspace of $\prod_\mathcal{U} L_1(\M)$. Since $L_1(\mathcal M_\mathcal U)$ is stable under $
\widetilde{\mathcal M}_\mathcal U$ actions, we know from \cite[Theorem III.2.7]{Tak79} that there must exist a central projection $e_\mathcal U$ in $\widetilde{\mathcal M}_\mathcal U$ such that $L_1(\mathcal M_\mathcal U) = L_1(\widetilde{\mathcal M}_\mathcal U )e_\mathcal U$. It turns out that $e_\mathcal U$ is the support projection of the trace $\tau_\mathcal U$. In the sequel, we will identify $\mathcal M_\mathcal U$ as a subalgebra of $\widetilde{\mathcal M}_\mathcal U$ by considering $\mathcal M_\mathcal U = \widetilde{\mathcal M}_\mathcal U e_\mathcal U$. More generally $$L_p(\mathcal M_\mathcal U) = L_p(\widetilde{\mathcal M}_\mathcal U)e_\mathcal U = \overline{\bigcup_{\widetilde{p}>p}L_{\widetilde{p}}(\widetilde{\mathcal M}_\mathcal U)}^{\|\cdot\|_{L_{{p}}(\widetilde{\mathcal M}_\mathcal U)}}$$ for all $0 < p < \infty$. Indeed, the subspace $L_p(\mathcal M_\mathcal U)$ can be characterized by using the notion of $p$-equiintegrability \cite{Ran02,Tak79}. The last identity above follows from this property, we refer to \cite[Lemma 1.7]{JuPe14} for further details. 

In this paper we shall make extensive use of column $L_p(\ell_2)$ spaces. In other words, the column subspaces of $S_p$-valued $L_p$-spaces. A more detailed definition can be found in \cite{Pis98}. The ultraproduct of column $L_p(\ell_2)$ spaces $$\prod_{\mathcal U}L_p(\mathcal M; \ell_2^c(\sigma))$$
forms a right $L_p$-module over $\widetilde{\mathcal M}_\mathcal U$ with right module action $$\xi\cdot x=(\xi^\sigma x^\sigma)^\bullet \quad \mbox{for} \quad x \in\prod_\mathcal U\mathcal M \quad \mbox{and} \quad \xi\in\prod_{\mathcal U}L_p(\mathcal M; \ell_2^c(\sigma)).$$
By \cite[Proposition 5.2]{JuSh05}, this module action extends naturally from $\prod_\mathcal U \mathcal M$ to $\widetilde{\mathcal M}_\mathcal U$ and does not depend on the representing families. Similarly, given $\xi = (\xi^\sigma)^\bullet$ and $\eta = (\eta^\sigma)^\bullet$ in $\prod_{\mathcal U}L_p(\mathcal M; \ell_2^c(\sigma))$, we consider the componentwise bracket $$\langle\xi,\eta\rangle_{\prod_{\mathcal U}L_p(\ell_2^c(\sigma))}= \big( \langle\xi^\sigma,\eta^\sigma\rangle \big)^\bullet= \Big( \sum_{t\in\sigma}\xi^{\sigma*}_t\eta^\sigma_t \Big)^\bullet\in \prod_\mathcal UL_{p/2}(\mathcal M) \simeq L_{p/2}(\widetilde{\mathcal M}_\mathcal U)$$
where $\xi^\sigma = \sum_{t\in\sigma} \xi^\sigma_t\otimes e_{t,0},\; \eta^\sigma = \sum_{t\in\sigma} \eta^\sigma_t\otimes e_{t,0}$. This defines an $L_{p/2}(\widetilde{\mathcal M}_\mathcal U)$-valued inner product which generates the norm of $\prod_{\mathcal U}L_p(\mathcal M; \ell_2^c(\sigma))$ and is compatible with the module action.
Hence $\prod_{\mathcal U}L_p(\mathcal M; \ell_2^c(\sigma))$ is a right $L_p$-module for $1\leq p\leq\infty$. The regular version of this right $L_p$-module is
$$X_p(\mathcal M) = \Big( \prod_{\mathcal U}L_p(\mathcal M; \ell_2^c(\sigma)) \Big) e_{\mathcal U}.$$
It also forms a right $L_p$-module over $\mathcal M_\mathcal U$ since for $\xi\in\prod_{\mathcal U}L_p(\mathcal M; \ell_2^c(\sigma))$, we have $\xi\in X_p(\mathcal M)$ iff $\langle\xi,\xi\rangle_{\prod_{\mathcal U}L_p(\ell_2^c(\sigma))} \in L_{p/2}({\mathcal M}_\mathcal U).$ Moreover, the family
$X_p(\mathcal M)_{1\leq p\leq\infty}$ forms a projective system of $L_p$-modules over $\mathcal M_\mathcal U$, see the arguments after the proof of Lemma 2.21 in \cite{JuPe14}.
As for the regular $L_p(\mathcal M_\mathcal U)$, we also have the following characterization 
\begin{align}\label{regular of lp module}
X_p(\mathcal M)=\overline{\bigcup_{\tilde{p}>p}\prod_{\mathcal U}L_{\tilde{p}}(\mathcal M; \ell_2^c(\sigma))}^{\prod_{\mathcal U}L_{{p}}(\ell_2^c(\sigma))}.
\end{align}
In the following, we shall also need the following subspaces
\begin{eqnarray*}
X^\circ_p(\M) & = & \Big\{(x^\sigma)^\bullet\in X_p(\M) \, \big| \; \mathcal{E}_t(x^\sigma_t)=0 \, \mbox{for all } t\in\sigma \Big\}, \\
X^{ad}_p(\M) & = & \Big\{(x^\sigma)^\bullet\in X_p(\M) \, \big| \; x^\sigma_t\in L_p(\mathcal{M}_t) \, \mbox{for all } t\in\sigma \Big\}.
\end{eqnarray*}

\begin{remark}
\emph{In what follows and in order to deal with diagonal Hardy spaces, we shall actually use the ultraproduct of the spaces $L_p(\mathcal M; \ell_2^c(\sigma\times\mathbb N))$ and its regularized version. We shall keep the (lighter) notation $X_p(\mathcal M)$ for this regularized version.}
\end{remark}

\subsection{Proof of Theorem Aii} This part is much more involved. We can not prove the desired estimate as a consequence of its validity in the discrete case. Instead we adapt the proof in the discrete case to give a direct argument. This roughly means that, for $x$ in the column Hardy spaces |algebraic atomic forms of the conditional and diagonal ones| we shall write $\E_t(x)$ in the form $A_tB$ with $A_t$ and $B$ in some amplified matrix algebra. Then we factorize $A_t$ using the symmetric Doob maximal inequality. In both steps of the argument, new difficulties appear which force us to establish some preliminary results in Lemmas \ref{lem:davis dec using hat} and \ref{lem:construction of conditional expectations} below.

\begin{lemma}\label{lem:davis dec using hat}
Let $1\leq p<2$ and $s \ge 2$ be determined by $1/p=1/2+1/s$. Any $x\in\M$ can be written in the form $x=x_c+x_{1_c}$, where the elements $x_c$ and $x_{1_c}$ satisfy in turn the following properties$\hskip1pt :$
\begin{itemize}
\item [i)] $x_c=w-L_p-\lim_{\sigma,\mathcal{U}}a^\sigma b^\sigma$ where
$$a^\sigma=\sum_{t\in\sigma}a^\sigma_t\otimes e_{1,t} \quad \mathrm{and} \quad b^\sigma=\sum_{t\in\sigma}b^\sigma_t\otimes e_{t,1}$$
satisfy the following estimate $\|(a^{\sigma*})^\bullet\|_{X^\circ_2}\|(b^\sigma)^\bullet\|_{X^{ad}_s}\leq C_{{p}}\|x\|_{\mathcal{H}^c_p}$.

\vskip5pt

\item[ii)] $x_{1_c}=w-L_p-\lim_{\sigma,\mathcal{U}} \big( \alpha^\sigma(1)\beta^\sigma(1)-\alpha^\sigma(2)\beta^\sigma(2) \big)$ where
$$\begin{array}{lcl} \hskip20pt \alpha^\sigma(1) = \displaystyle \sum_{t\in\sigma}u_t(\alpha^{\sigma*}_t)^*\otimes e_{1,t}, & & \beta^\sigma(1) = \displaystyle  \sum_{t\in\sigma}u_t(\beta^\sigma_t)\otimes e_{t,1}, \\ \hskip20pt \alpha^\sigma(2) = \displaystyle \sum_{t\in\sigma}u_{t^-(\sigma)}(\alpha^{\sigma*}_t)^*\otimes e_{1,t}, & & \beta^\sigma(2) = \displaystyle  \sum_{t\in\sigma}u_{t^-(\sigma)}(\beta^\sigma_t)\otimes e_{t,1}, \end{array}$$ for some module maps $u_t$ and $u_{t^-(\sigma)}$ acting on $\alpha_t^{\sigma*}, \beta_t^\sigma$. We also have
\begin{eqnarray*}
\hskip20pt \lefteqn{\max \Big\{ \big\| (\alpha^\sigma(1)^*)^\bullet \big\|_{X_2} \big\| (\beta^\sigma(1))^\bullet \big\|_{X_s},
\big\| (\alpha^\sigma(2)^*)^\bullet \big\|_{X_2} \big\| (\beta^\sigma(2))^\bullet \big\|_{X_s} \Big\}} \\
& \hskip20pt \le & C_s \lim_{\sigma,\mathcal{U}} \Big\| \sum_{t\in\sigma}\alpha^{\sigma*}_t \otimes e_{t,1} \Big\|_{2}
\, \lim_{\sigma,\mathcal{U}} \Big\| \sum_{t\in\sigma} \beta^\sigma_t \otimes e_{t,1} \Big\|_s \ \le \ C_{p,s} \|x\|_{\mathcal{H}^c_p}.
\end{eqnarray*}
\end{itemize}
\end{lemma}

\dem
Assume by homogeneity that $\|x\|_{\mathcal{H}^c_p}<1$. According to \cite[Lemma 3.2.25]{Per} we can find $p < \widetilde{p} < 2$ such that $\|x\|_{\mathcal{H}^c_{\widetilde{p}}}<1$. Applying Davis decomposition \cite{JM1,Per0} to every partition $\sigma$, we find a decomposition  $x=x_{c}^\sigma+x_{1_c}^\sigma$ such that
$$\|x_c^\sigma\|_{h^c_{\widetilde{p}}(\sigma)}+\|x_{1_c}^\sigma\|_{h^{1_c}_{\widetilde{p}}(\sigma)}<C_{\widetilde{p}}$$
where $C_{\widetilde{p}}$ is the optimal constant in the Davis decomposition.
It is easy to see that $x_c^\sigma$ and $x_{1_c}^\sigma$ are uniformly bounded in $L_{\widetilde{p}}(\mathcal M)$ since $h^c_{\widetilde{p}}(\sigma)$ and $h^{1_c}_{\widetilde{p}}(\sigma)$ embed into $L_{\widetilde{p}}(\mathcal M)$ uniformly. By \cite[Lemma 3.1.6]{Per}, we can define $$x_c=w-L_p-\lim_{\sigma,\mathcal{U}}x_c^{\sigma}\quad \mathrm{and} \quad x_{1_c}=w-L_p-\lim_{\sigma,\mathcal{U}}x_{1_c}^{\sigma}.$$ Let us check the properties of $x_c$. According to the algebraic atomic decomposition in the discrete case |see for instance \cite[Lemma 2.1]{HJP}| for any $\delta>0$ we find a decomposition $$x_c^\sigma=\sum_{t\in\sigma}a^\sigma_tb^\sigma_t = \Big( \sum_{t\in\sigma}a^\sigma_t\otimes e_{1,t} \Big) \Big( \sum_{t\in\sigma}b^\sigma_t\otimes e_{t,1} \Big) = a^\sigma b^\sigma$$ such that $\mathcal{E}_ta^\sigma_t=0$, $b^\sigma_t\in L_{\widetilde{s}}(\mathcal{M}_t)$ where $1/{\widetilde{p}}=1/2+1/{\widetilde{s}}$ and
$$\|a^\sigma\|_2 \|b^\sigma\|_{\widetilde{s}}\leq (1+\delta)\|x_c^\sigma\|_{h^c_{\widetilde{p}}(\sigma)}.$$ By \eqref{regular of lp module} and since $b^\sigma$ are uniformly bounded in $L_{\widetilde{s}}(\mathcal M;\ell^c_2({\sigma}))$, we conclude that $b=(b^\sigma)^\bullet$ belongs to $X^{ad}_{s}(\mathcal M)$. This implies that  $b=be_{\mathcal{U}}$ and recalling that $e_{\mathcal{U}}$ is a central projection we get $(x_c^\sigma)^\bullet=(a^\sigma (b^\sigma e^\sigma_{\mathcal{U}}))^\bullet=\big((e^\sigma_{\mathcal{U}} a^\sigma) (b^\sigma e^\sigma_{\mathcal{U}})\big)^\bullet$. Therefore if we reset $a^*$ to be $a^*e_{\mathcal{U}}$, we obtain the desired properties of $x_c$.

The properties of $x_{1_c}$ can be checked similarly, we shall only give a sketchy argument. According to the definition of $h^{1_c}_{\widetilde{p}}(\sigma)$ |see \cite[Lemma 2.1]{HJP}| for any $\delta>0$, we have $x_{1_c}^\sigma=\sum_{t\in\sigma}d^\sigma_t(\alpha^\sigma_t\beta^\sigma_t)$ with $d^\sigma_t=\mathcal E_t-\mathcal E_{t^-(\sigma)}$ the martingale difference associated to $\sigma$ and $$\Big\| \sum_{t\in\sigma}\alpha^\sigma_t\otimes e_{1,t} \Big\|_2 \Big\| \sum_{t\in\sigma}\beta^\sigma_t\otimes e_{t,1} \Big\|_{\widetilde{s}} \le (1+\delta)\|x_{1_c}^\sigma\|_{h^{1_c}_{\widetilde{p}}(\sigma)}$$ where $1/\widetilde{p}=1/2+1/\widetilde{s}$. Now applying Proposition 2.8 of \cite{Jun02}, we write 
\begin{eqnarray*}
\mathcal E_t(\alpha^\sigma_t\beta^\sigma_t) & = & u_t(\alpha^{\sigma*}_t)^*u_t(\beta^\sigma_t), \\ \mathcal E_{t^-(\sigma)}(\alpha^\sigma_t\beta^\sigma_t) & = & u_{t^-(\sigma)}(\alpha^{\sigma*}_t)^*u_{t^-(\sigma)}(\beta^\sigma_t),
\end{eqnarray*}
where $u_s:\mathcal M\rightarrow C(\mathcal M_s)$ is an isometric right $\mathcal M_s$-module map with $s=t,t^-(\sigma)$.
Then we further write $$x_{1_c}^\sigma=\sum_{t\in\sigma}u_t(\alpha^{\sigma*}_t)^*u_t(\beta^\sigma_t) - \sum_{t\in\sigma} u_{t^-(\sigma)}(\alpha^{\sigma*}_t)^*u_{t^-(\sigma)}(\beta^\sigma_t).$$ The rest of the argument is similar to the one used for $x_c$. The only significant difference is that the inequality requires the dual form of Doob's inequality. \fin

To prove the asymmetric maximal inequality for conditional Hardy spaces, we need the following lemma, which provides an explicit representation of the sequence of right $L_p$-modules $(X_p(\M_t))_{t\in\sigma}$ and the resulting module maps extend naturally to conditional expectations associated to certain von Neumann algebras. 

\begin{lemma}\label{lem:construction of conditional expectations}
Given $0\leq s<t\leq1$, we can find Hilbert spaces $H_s\subset H_t$, two projections $e_s\in \mathcal{B}(H_s)\bar{\otimes}\M_s$, $e_t\in\mathcal{B}(H_t)\bar{\otimes}\mathcal{M}_t$, and one module homomorphism $\iota$ satisfying 
\begin{eqnarray*}
\iota(X_p(\M_s)) & = & e_s L_p(\mathcal{B}(H_s)\bar{\otimes}\M_s) e_{1,1}, \\ 
\iota(X_p(\M_t \hskip1pt )) & = & e_t  \hskip1pt L_p(\mathcal{B}(H_t)\bar{\otimes}\M_t \hskip1pt ) e_{1,1},
\end{eqnarray*}
for $1\leq p\leq\infty$. Furthermore, we find a $\mathcal M_s$-module map $$\mathcal{P}_s: e_t  L_p(\mathcal{B}(H_t)\bar{\otimes}\M_t) e_{1,1}\rightarrow e_s L_p(\mathcal{B}(H_s)\bar{\otimes}\M_s) e_{1,1}$$ which extends to a conditional conditional expectation given by $$\widehat{\mathbb{E}}_s: e_t\mathcal{B}(H_t)\bar{\otimes}\M_t e_t \to e_t \big( \mathcal{B}(H_s)\bar{\otimes}\M_s\oplus \mathcal{D}_{t-s}\bar{\otimes}\M_t \big) e_t$$ where $\mathcal{D}_{t-s}$ is the diagonal subalgebra of the matrix algebra $\mathcal{B}(H_t\ominus H_s)$. In fact, by induction, the same construction applies for all elements at once of any given partition $\sigma \in \Sigma$. This yields the following maps for all $t \in \sigma$ $$\iota^\sigma: X_p(\M_t) \to L_p(\mathcal{B}(H_t)\bar{\otimes}\M_t ) e_{1,1},$$ $$\mathcal{P}_t: e_1  L_p(\mathcal{B}(H_1)\bar{\otimes}\M) e_{1,1}\rightarrow e_t L_p(\mathcal{B}(H_t)\bar{\otimes}\M_t) e_{1,1},$$ $$\widehat{\mathbb{E}}_t: e_1 \mathcal{B}(H_1)\bar{\otimes}\M e_1 \to e_1 \big( \mathcal{B}(H_t)\bar{\otimes}\M_t\oplus \mathcal{D}_{1-t}\bar{\otimes}\M \big) e_1.$$
\end{lemma}

\dem According to \cite[Theorem 2.5]{JuSh05}, any $L_p$-module is isometrically isomorphic to a principal $L_p$-module. In particular, there exists a sequence of projections $(q_\alpha)_{\alpha\in{I_s}}$ in $\mathcal M_s$ so that $X_p(\M_s)$ is isometrically isomorphic to $\bigoplus q_\alpha L_p(\M_s)$ as an $L_p$-module. Moreover, we can assume from the proof in \cite{JuSh05} that $q_1=\1$, the identity in $\mathcal M_s$. Hence, we can formally write $$x\sim \bigoplus q_\alpha u_\alpha(x).$$ On the other hand, since $(X_p(\M_s))_{1\leq p\leq\infty}$ is a projective system of $\mathcal{M}_s$-modules \cite[Lemma 2.21]{JuPe14}, $(q_\alpha)_{\alpha \in {I_s}}$ can be chosen to be independent of $p$. Consider the Hilbert space $H_s = \ell_2(I_s)$ and define $e_s = \sum_{\alpha\in I_s}e_{\alpha,\alpha}\otimes q_\alpha$. Then, the module homomorphism $\iota$ is given by $$\iota: X_p(\mathcal{M}_s) \ni x \mapsto \sum_{\alpha\in I_s}e_{\alpha,1}\otimes q_\alpha u_{\alpha}(x) \in e_s L_p(\mathcal{B}(H_s)\bar{\otimes}\M_s) e_{1,1}.$$ Since $\mathcal{M}_s$ is a von Neumann subalgebra of $\M_t$ for  $s<t$, we complete the basis $(q_\alpha)_{\alpha\in I_s}$ in $\mathcal{M}_t$ as in the proof of \cite[Theorem 2.5]{JuSh05} to obtain a basis $(q_\alpha)_{\alpha\in I_t}$ in $\M_t$ such that $X_p(\M_t)$ is isometrically isomorphic to $\bigoplus q_\alpha L_p(\M_t)$. Therefore, we may reproduce the construction above with $H_t \supset H_s$ and $e_t \ge e_s$. The module map $\mathcal{P}_s$ is then constructed as follows $$\mathcal{P}_s \Big( \sum_{\alpha\in I_t}e_{\alpha,1}\otimes q_\alpha m_\alpha \Big) = \sum_{\alpha\in I_s} e_{\alpha,1} \otimes q_\alpha \mathcal{E}_s (m_\alpha).$$ The related conditional expectation $\widehat{\mathbb{E}}_s$ is simply defined as
\begin{eqnarray*}
\widehat{\mathbb{E}}_s \big( e_t (m_{\alpha,\beta})_{\alpha,\beta\in I_t} e_t \big) & = & \widehat{\mathbb{E}}_s \big( (q_\alpha m_{\alpha,\beta} q_\beta)_{\alpha,\beta\in I_t} \big) \\ & = & e_s(\mathcal{E}_s(m_{\alpha,\beta}))_{\alpha,\beta\in I_s}e_s\oplus\sum_{\alpha\in I_t\setminus I_s} e_{\alpha,\alpha} \otimes q_\alpha m_{\alpha,\alpha}q_{\alpha}.
\end{eqnarray*}
Since $q_1=\1$, we observe that $$e_{1,1} \otimes \1 \le e_s.$$ Therefore, it is easy to check that $\widehat{\mathbb{E}}_s(x) = \mathcal{P}_s(x)$ whenever $x \in e_t \mathcal{B}(H_t)\bar{\otimes}\M_t e_{1,1}$. The last assertion follows inductively just replacing $(s,t)$ by $(t,1)$ for all $t \in \sigma$. \fin

\noindent \textbf{Proof of Theorem Aii.} By Lemma \ref{lem:davis dec using hat}, it suffices to prove
\begin{eqnarray}
\label{estimate hhcp}
\|(\E_t(x_c))_t\|_{L_p(\ell^\theta_\infty)} & \le & C(p,\theta)\|x\|_{\mathcal{H}^c_p}, \\
\label{estimate hh1cp}
\|(\E_t(x_{1_c}))_t\|_{L_p(\ell^\theta_\infty)} & \le & C(p,\theta)\|x\|_{\mathcal{H}^c_p}.
\end{eqnarray}
Given a fixed $t \in [0,1]$ and by the definition of ultralimit, we may clearly restrict our ultralimits to run over partitions $\sigma$ satisfying $t \in \sigma$. Moreover, the self-adjointness and contractivity of conditional expectations allow us to commute them with the weak $L_p$ ultralimits below. Finally, emulating our argument in the discrete case \cite{HJP} we may combine the mean zero property of $X^\circ_2(\mathcal M)$ with the adapted sequences in $X^{ad}_{s}(\mathcal M)$ to conclude for $\widehat{\mathcal{E}}_t=(\mathcal{E}_t\otimes id_{\mathcal{B}(\ell_2(\sigma))})^{\bullet}$ that 
\begin{eqnarray*}
\mathcal{E}_t(x_c) & = & \mathcal{E}_t \Big( w-L_p-\lim_{\sigma,\mathcal{U}} \sum_{s \in\sigma} a^\sigma_s b^\sigma_s \Big) \ = \ w-L_p-\lim_{t\in\sigma,\mathcal{U}}\mathcal{E}_t \Big( \sum_{s\in\sigma}a^\sigma_sb^\sigma_s \Big) \\ & = & w-L_p-\lim_{t\in\sigma,\mathcal{U}}\widehat{\mathcal{E}}_t(a^\sigma)b^\sigma \ = \ w-L_p-\lim_{\sigma,\mathcal{U}}\widehat{\mathcal{E}}_t(a^\sigma)b^\sigma \ = \ \big\langle \widehat{\mathcal{E}}_t(a^*),b \big\rangle
\end{eqnarray*}
for the right module bracket $\langle \cdot , \cdot \rangle$. Using the argument in the proof of Lemma \ref{lem:construction of conditional expectations}, there must exist a Hilbert space $K$, a projection $f\in \mathcal{B}(K)\bar{\otimes}\M$ and a module homomorphism $\iota$ such that for all $1\leq q\leq\infty$ $$\iota(X_q(\M)) = f L_p(\mathcal{B}(K) \bar{\otimes} \M) e_{1,1}.$$ Therefore, it turns out that for any $t\in[0,1]$ we have $$\mathcal{E}_t(x_c) = \iota(\widehat{\mathcal{E}}_t(a^*))^* \iota(b).$$ The desired factorization of $\mathcal{E}_t(x_c)$ shall be deduced from a factorization of $\iota(\widehat{\mathcal{E}}_t(a^*))$ in the von Neumann algebra $f(\mathcal{B}(K)\bar{\otimes}\M)f$ which is based on the following estimate
\begin{equation} \label{1 est in h2pc}
\big\|(\iota(\widehat{\mathcal{E}}_t(a^*)))_{t\in[0,1]}\big\|_{L_2(f(\mathcal{B}(K)\bar{\otimes}\M)f;\ell^{\rho}_\infty)}\leq C_{\rho}\|a^*\|_{X^\circ_2}
\end{equation}
where $0 < \rho < 1$ is the real number satisfying $2/\rho=p/(1-\theta)$. Indeed, let us complete the proof of the inequality for $x_c$ assuming that \eqref{1 est in h2pc} holds. This estimate yields a factorization $\iota(\widehat{\mathcal{E}}_t(a^*)) = yv_tz$ in $$L_{\frac{2}{1-\rho}}(\mathcal{B}(K)\bar{\otimes}\M)\times L_{\infty}(\mathcal{B}(K)\bar{\otimes}\M)\times L_{\frac{2}{\rho}}(\mathcal{B}(K)\bar{\otimes}\M),$$ where the $v_t$'s are contractions in $\mathcal{B}(K)\bar{\otimes}\M$ and $$\|y\|_{\frac{2}{1-\rho}} \|z\|_{\frac{2}{\rho}} \le C_{\rho} \|a^*\|_{X^\circ_2}.$$ The promised factorization of $\mathcal{E}_t(x_c)$ is then given by $$\mathcal{E}_t(x_c)=\underbrace{(z^*z)^{1/2}}_{\alpha}\underbrace{(z^*z)^{-1/2}z^*v_t^*y^*\iota(b)(\iota(b)^*yy^*\iota(b))^{-1/2}}_{w_t} \underbrace{(\iota(b)^*yy^*\iota(b))^{1/2}}_{\beta}.$$ Namely, note that that image of $\iota$ is composed of column matrices, so the same holds for $z$ and $\iota(b)$. In particular, $\alpha, \beta$ and $w_t$ are affiliated to $\M \subsetneq \mathcal{B}(K) \bar\otimes \M$ and it suffices to prove that $(\alpha, w_t, \beta) \in L_{p/(1-\theta)}(\M) \times L_\infty(\M) \times L_{p/\theta}(\M)$ with product norm dominated by the $\mathcal{H}_p^c$-norm of $x$. This easily follows from the estimates above, H\"older's inequality and Lemma \ref{lem:davis dec using hat} $$\|\alpha\|_{\frac{p}{1-\theta}} \Big( \sup_{t\in[0,1]} \|w_t\|_\infty \Big) \|\beta\|_{\frac{p}{\theta}} \le c_{p,\theta}\|a^*\|_{X^\circ_2}\|b\|_{X^{ad}_s} \le C(p,\theta) \|x\|_{\mathcal{H}_p^c}.$$
Theorefore, \eqref{estimate hhcp} will follow if we justify our estimate \eqref{1 est in h2pc}. Note that $L_2(\mathcal M;\ell^{\rho}_\infty)$ is a Banach space, using similar arguments as those used in \cite[Proposition 2.1]{JuXu07} we get $$\big\| (\iota(\widehat{\mathcal{E}_t}(a^*)))_{t\in[0,1]} \big\|_{L_2(\ell^{\rho}_\infty)} = 
\sup_{\sigma \in \Sigma} \big\| (\iota(\widehat{\mathcal{E}_t}(a^*)))_{t \in\sigma} \big\|_{L_2(\ell^{\rho}_\infty)}.$$ Given a partition $\sigma$, we know from \cite[Corollary 5.2]{JuSh05} that  there exists a partial isometry $u^\sigma$ satisfying $\iota=u^\sigma\iota^\sigma$ for the module homomorphism $\iota^\sigma$ constructed in Lemma \ref{lem:construction of conditional expectations} above. Moreover
\begin{eqnarray*}
\big\| (u^\sigma\iota^\sigma(\widehat{\mathcal{E}_t}(a^*)))_{t\in\sigma} \big\|_{L_2(f(\mathcal{B}(K)\bar{\otimes}\M)f;\ell^{\rho}_\infty)} \!\!\! & \le & \!\!\! \big\| (\iota^\sigma(\widehat{\mathcal{E}_t}(a^*)))_{t\in\sigma} \big\|_{L_2(e_1(\mathcal{B}(H_1)\bar{\otimes}\M)e_1;\ell^{\rho}_\infty)}\\
\!\!\! & = & \!\!\! \big\| (\widehat{\mathbb{E}}_t(\iota^\sigma(a^*)))_{t\in\sigma} \big\|_{L_2(e_1(\mathcal{B}(H_1)\bar{\otimes}\M)e_1;\ell^{\rho}_\infty)},
\end{eqnarray*}
where we have used above the intertwining identity $\iota^\sigma\circ\widehat{\mathcal{E}_t}=\widehat{\mathbb{E}}_t \circ\iota^\sigma$ which follows as in \cite[Theorem 1.5]{JuSh05}, see also \cite{Pas73}. Finally we may apply the asymmetric Doob maximal inequalities in discrete time \cite[Theorem A]{HJP} for the increasing sequence of conditional expectations from Lemma \ref{lem:construction of conditional expectations}. This yields $$\big\| (\iota(\widehat{\mathcal{E}}_t(a^*)))_{t\in[0,1]} \big\|_{L_2(\ell^{\rho}_\infty)} \le C_\rho \sup_{\sigma \in \Sigma} \big\| \iota^\sigma(a^*) \big\|_{2} = C_{\rho} \|a^*\|_{X^\circ_2}.$$

\noindent In order to prove \eqref{estimate hh1cp} and arguing as above, we start by noticing that
\begin{eqnarray*}
\mathcal{E}_t(x_{1_c}) & = & w-L_p-\lim_{t\in\sigma,\mathcal{U}}\sum_{s\in\sigma,s\leq t}d^\sigma_s(\alpha^\sigma_s\beta^\sigma_s) \\
& = & w-L_p-\lim_{t\in\sigma,\mathcal{U}}\sum_{s\in\sigma,s\leq t}\mathcal{E}_s(\alpha^\sigma_s\beta^\sigma_s) - \sum_{s\in\sigma,s\leq t}\mathcal{E}_{s^-(\sigma)}(\alpha^\sigma_s\beta^\sigma_s) \ = \ Y_t-Z_t
\end{eqnarray*}
from (the proof of) Lemma \ref{lem:davis dec using hat}. By the triangle inequality, we shall deal with $Y_t$ and $Z_t$ separately. Since both can be handled similarly, we restrict our attention to $Y_t$. Define the conditional expectations $$\mathbf{E}^\sigma_t=id_{\mathcal{M}}\otimes id_{\mathcal{B}(\ell_2)}\otimes \mathsf{E}^\sigma_t,$$ where $$\mathsf{E}^\sigma_t \big( (m_{u,v})_{u,v\in\sigma} \big) = \big( m_{u,v} \big)_{u,v\leq t} \oplus \sum_{u > t} m_{u,u} \otimes e_{u,u}.$$ Then, it is easily checked that we can write
\begin{eqnarray*}
Y_t & = & w-L_p-\lim_{t\in\sigma,\mathcal{U}}\mathbf{E}^\sigma_t(\alpha^\sigma(1))\beta^\sigma(1) \\
& = & w-L_p-\lim_{\sigma,\mathcal{U}}\mathbf{E}^\sigma_t(\alpha^\sigma(1))\beta^\sigma(1) \ = \ \big\langle \widehat{\mathbf{E}}_t(\alpha(1)^*),\beta(1) \big\rangle
\end{eqnarray*}
where $\widehat{\mathbf{E}}_t = (\mathbf{E}^\sigma_t)^\bullet$ and $\langle \cdot, \cdot \rangle$ stands for the right module bracket. As proved in Lemma \ref{lem:construction of conditional expectations}, for any fixed partition $\sigma \in \Sigma$ we can construct a Hilbert space $H$, a projection $e$ and a module homomorphism $\iota^\sigma$ such that $$\iota^\sigma(X_q(\mathcal{M}))=eL_q(\mathcal{B}(H)\bar{\otimes}\M)e_{1,1} \quad \mbox{for any} \quad 1\leq q\leq\infty.$$ On the other hand, it is easy to check that $\widehat{\mathbf{E}}_t$ is a right module map. Then by the characterizations of module maps \cite{JuSh05,Pas73},  there exists a projection $p_t\in \mathcal{B}(H)\bar{\otimes}\M$ such that $e_{1,1}\otimes 1\leq p_t\leq e$ and $$\iota^\sigma(\widehat{\mathbf{E}}_t(z))=p_t\iota^\sigma(z)\in p_tL_q(\mathcal{B}(H)\bar{\otimes}\M)e_{1,1} \quad \mbox{for all} \quad z \in X_q(\M).$$ Consider the increasing von Neumann subalgebras  
$$A_t = p_t(\mathcal{B}(H)\bar{\otimes}\M)p_t\oplus  p^\perp_t(\mathcal{D}(H)\bar{\otimes}\M)p^\perp_t \subset e(\mathcal{B}(H)\bar{\otimes}\M)e$$ for $t\in\sigma$ and $\sigma$ fixed, where $\mathcal{D}(H)$ is the diagonal subalgebra of $\mathcal{B}(H)$. The associated conditional expectations are $\mathbb{E}_t(z)=p_tzp_t\oplus p^\perp_t\mathcal{D}(z)p^\perp_t$ for every element $z \in e(\mathcal{B}(H)\bar{\otimes}\M)e$ with diagonal part $\mathcal{D}(z)$. Then it is easy to check that when $z \in eL_q(\mathcal{B}(H)\bar{\otimes}\M)e_{1,1}$, we have $$p_t\iota^\sigma(z)=\mathbb{E}_t(\iota^\sigma(z)).$$ This enables us to follow the arguments in the proof of estimate \eqref{estimate hhcp}. \fin

\section{{\bf Proof of Theorem B}} 

Motivated by our proof of Theorem Aii and particularly Lemma \ref{lem:davis dec using hat}, we first introduce new families of Hardy spaces. Let $1\leq p<2$, and $w,s\geq2$ such that $1/p=1/w+1/s$. We define
\begin{eqnarray*}
\hskip20pt \widehat{h}^{c}_{pw}(\mathcal M) & = & \big\{ x \in L_p(\M) : \|x\|_{\widehat{h}^c_{pw}}<\infty \big\}, \\
\widehat{h}^{1_c}_{pw}(\mathcal M) & = & \big\{ x \in L_p(\M) : \|x\|_{\widehat{h}^{1_c}_{pw}}<\infty \big\},
\end{eqnarray*}
with
\begin{eqnarray*}
\|x\|_{\widehat{h}^{c}_{pw}} & = & \inf_{x=\langle a^*,b\rangle}\|a^*\|_{X^\circ_w}\|b\|_{X^{ad}_s}, \\
\|x\|_{\widehat{h}^{1_c}_{pw}} & = & \inf_{x=\langle\langle a^*,b\rangle\rangle}\|a^*\|_{X_w}\|b\|_{X_s},
\end{eqnarray*}
where the inner products are defined as
\begin{eqnarray*}
\hskip20pt \langle a^*,b\rangle & = & w-L_p-\lim_{\sigma,\mathcal{U}}\sum_{t\in\sigma}a^\sigma_{t}b^\sigma_{t}, \\
\langle\langle a^*,b\rangle\rangle & = & w-L_p-\lim_{\sigma,\mathcal{U}}\sum_{t\in\sigma}d^\sigma_{t}(a^\sigma_{t}b^\sigma_{t}).
\end{eqnarray*}
The row spaces $\widehat{h}^{r}_{pw}(\mathcal M)$ and $\widehat{h}^{1_r}_{pw}(\mathcal M)$ are defined in a similar way.

\begin{remark}\label{norm of hat hardy}
\emph{Both families of Hardy spaces are Banach spaces. The properties of the norm follow as in the discrete case \cite[Lemma 2.3]{HJP}, while the completeness requires some details from the proof of Theorem \ref{pro:hat hcpw=hcpw}, see also \cite[Lemma 2.1]{HJP}.}
\end{remark}

\subsection{Proof of Theorem Bi} The inclusion $$\widehat{\H}_{pw}^c(\M) = \widehat{h}^{c}_{pw}(\M)+\widehat{h}^{1_c}_{pw}(\M)\subset L_p(\M)$$ trivially holds by definition and the fact that it was already proved in the discrete case \cite[Lemma 2.5]{HJP}. Let us prove the reverse inclusion. Given $x\in \M$ with $\|x\|_p<1$, there exists $p<\widetilde{p}<2$ such that $\|x\|_{\widetilde{p}}<1$. Applying \cite[Theorem Bi]{HJP} to $\widetilde{p}$ and each partition $\sigma$, we get a decomposition $x=x_c^{\sigma}+x_{1_c}^{\sigma}+x_r^{\sigma}+x_{1_r}^{\sigma}$ with $x_u^\sigma\in h^u_{{\widetilde{p}}\widetilde{w}}(\sigma)$ for $u=c,1_c,r,1_r$ and certain $\widetilde{w}>w$, $\widetilde{s}>s$ so that 
$$\|x_c^{\sigma}\|_{h^c_{\widetilde{p}\widetilde{w}}(\sigma)}+\|x_{1_c}^{\sigma}\|_{h^{1_c}_{\widetilde{p}\widetilde{w}}(\sigma)}
+\|x_r^{\sigma}\|_{h^{r}_{\widetilde{p}\widetilde{w}}(\sigma)}+\|x_{1_r}^{\sigma}\|_{h^{1_r}_{\widetilde{p}\widetilde{w}}(\sigma)} \le C_{\widetilde{p},\widetilde{w}} \|x\|_{\widetilde{p}} < C_{\widetilde{p},\widetilde{w}}.$$ Moreover, according to \cite[Lemma 2.5]{HJP} we conclude that the family $(x_u^\sigma)_{\sigma \in \Sigma}$ is uniformly bounded in the reflexive space $L_{\widetilde{p}}(\M)$; whence in $L_p(\M)$ and we can define $$x_u=w-L_p-\lim_{\sigma,\mathcal{U}}x_u^\sigma\in L_{{p}}.$$ Now, since $\widetilde{w}>w$ and $\widetilde{s}>s$, we may use \eqref{regular of lp module} and pick $$a^*\in X^\circ_w(\mathcal M) \quad \mbox{and} \quad b\in X^{ad}_s(\mathcal M)$$ such that $x_c=\langle a^*,b\rangle$ and $$\|a^*\|_{X^\circ_w}\|b\|_{X^{ad}_s}\leq 2\lim_{\sigma,\mathcal{U}}\|x_c^\sigma\|_{h^c_{pw}(\sigma)}\leq 2\lim_{\sigma,\mathcal{U}}\|x_c^\sigma\|_{h^c_{\widetilde{p}\widetilde{w}}(\sigma)}\leq2C_{\widetilde{p},\widetilde{w}}.$$
Therefore we get $x_c\in\widehat{h}^c_{pw}(\mathcal M)$. Similar arguments also work for $u=1_c,r,1_r$. \fin

\subsection{Proof of Theorem Bii} It is straightforward to adapt  the proof of Theorem Aii to the present setting by using the new Hardy spaces. The only significant difference is that we assume $w>2$ to obtain the completely asymmetric estimates. Instead of \eqref{1 est in h2pc}, we use the following estimate from \cite[Corollary 4.6]{Jun02} $$\big\| (\iota(\widehat{\mathcal{E}}_t(a^*)))_{t\in[0,1]} \big\|_{L_2(f(\mathcal{B}(K)\bar{\otimes}\M)f;\ell^{r}_\infty)} \le C \|a^*\|_{X^\circ_2}.$$

\subsection{Proof of Theorem Biii} Unlike the discrete case, the proof of this part becomes quite involved since it requires an algebraic atomic characterization of the hat Hardy spaces, which could be regarded as a continuous version of Lemma 2.1 in \cite{HJP}. We shall need some preliminaries. Let $1\leq p<2$ and $w,s\geq2$ such that $1/p=1/w+1/s$. Given $x\in\M$, we consider the following quantities and their corresponding row analogues $$\|x\|_{h^c_{pw}}=\lim_{\sigma,\mathcal{U}}\|x\|_{h^c_{pw}(\sigma)} \quad \mathrm{and} \quad \|x\|_{h^{1_c}_{pw}}=\lim_{\sigma,\mathcal{U}}\|x\|_{h^{1_c}_{pw}(\sigma)}.$$ These limits exists |see Lemma \ref{lem:monotonicity} below| and define four norms.

\begin{lemma}\label{lem:monotonicity}
Let $1\leq p<2$ with $$\frac1p = \frac1w + \frac1s \quad \mbox{for some} \quad w,s \ge 2.$$ Then, the following inequalities hold for $\sigma \subset \sigma'$ 
\begin{eqnarray*}
\|x\|_{h^c_{pw}(\sigma)} & \ge & \|x\|_{h^c_{pw}(\sigma')}, \\ 
\|x\|_{h^{1_c}_{pw}(\sigma)} & \le & C_{w,s} \|x\|_{h^{1_c}_{pw}(\sigma')}.
\end{eqnarray*}
Similar results hold in the row case. In particular, we get 
$$\begin{array}{rclcrcl} \|x\|_{h^r_{pw}} & = & \displaystyle \inf_{\sigma \in \Sigma}\|x\|_{h^r_{pw}(\sigma)}, & & \|x\|_{h^{1_r}_{pw}} & \sim & \displaystyle   \sup_{\sigma \in \Sigma} \|x\|_{h^{1_r}_{pw}(\sigma)}, \\ \|x\|_{h^c_{pw}} & = & \displaystyle  \inf_{\sigma \in \Sigma}\|x\|_{h^c_{pw}(\sigma)}, & & \|x\|_{h^{1_c}_{pw}} & \sim & \displaystyle   \sup_{\sigma \in \Sigma} \|x\|_{h^{1_c}_{pw}(\sigma)}. \end{array}$$
\end{lemma}

\dem The assertion for $h_{pw}^c$ follows trivially since any $x=\sum_{t\in\sigma} a_t^\sigma b_t^\sigma$ admits a representation in the form $x=\sum_{t\in\sigma'} a_t^\sigma b_t^\sigma$ by taking $a_t^\sigma = b_t^\sigma = 0$ for $t \notin \sigma$. The space $h_{pw}^{1_c}$ requires a more involved argument. Given $\varepsilon>0$, let $x \in \M$ admit a decomposition $$x = \sum_{s \in \sigma'} d^{\sigma'}_s(a_s^{\sigma'} b_s^{\sigma'})$$ satisfying the following estimate $$\Big\| \Big( \sum_{s \in \sigma'} |a_s^{\sigma'*}|^2 \Big)^{\frac12} \Big\|_w \Big\| \Big( \sum_{s \in\sigma'} |b_s^{\sigma'}|^2 \Big)^{\frac12} \Big\|_s \le \|x\|_{h^{1_c}_{pw}(\sigma')}+\varepsilon.$$ It is easy to check that we can rewrite $$x = \sum_{t\in\sigma} d^{\sigma}_t \Big( \sum_{s\in J_t} d^{\sigma'}_s(a_s^{\sigma'}b_s^{\sigma'}) \Big),$$ where $J_t$ denotes the collection of $s\in\sigma'$ such that $t^-(\sigma)\leq s^-(\sigma')<s\leq t$. By triangle inequality, $\|x\|_{h^{1_c}_{pw}(\sigma)}$ is smaller than the sum of the following two quantities
$$\Big\| \sum_{t\in\sigma} d^{\sigma}_t \Big( \sum_{s\in J_t} \mathcal{E}_s(a_s^{\sigma'} b_s^{\sigma'}) \Big) \Big\|_{h^{1_c}_{pw}(\sigma)} + \Big\| \sum_{t \in\sigma} d^{\sigma}_t \Big( \sum_{s \in J_t} \mathcal{E}_{s^-(\sigma')}(a_s^{\sigma'} b_s^{\sigma'}) \Big) \Big\|_{h^{1_c}_{pw}(\sigma)}.$$ We only estimate the first term, the second one follows similarly. Define $$a_t^\sigma = \sum_{s \in J_t} \mathcal{E}_s(a_s^{\sigma'} b_s^{\sigma'})(b_t^\sigma)^{-1} \quad \mbox{where} \quad b_t^\sigma = \Big( \sum_{s \in J_t} \mathcal{E}_s |b_s^{\sigma'}|^2 \Big)^{\frac12} \quad \mbox{for} \quad t\in\sigma.$$ Here we are assuming by approximation that $b_t^\sigma$ is invertible. We can rewrite $$x = \sum_{t\in\sigma} d^{\sigma}_t \Big( \sum_{s \in J_t} \mathcal{E}_s(a_s^{\sigma'} b_s^{\sigma'}) \Big) = \sum_{t\in\sigma} d^{\sigma}_t (a_t^\sigma b_t^\sigma).$$ Note that $1 \le p < 2$ and $w,s \ge 2$ by assumption, which implies $2 \le s < \infty$. In particular, the dual version of Doob's maximal inequality in $L_{s/2}(\M)$ yields the following estimate 
\begin{eqnarray*}
\big\| (b_t^\sigma)_{t\in\sigma} \big\|_{L_s(\ell^c_2)} & = & \Big\| \Big( \sum_{t \in\sigma} \sum_{s \in J_t} \mathcal{E}_s |b_s^{\sigma'}|^2 \Big)^{\frac12} \Big\|_s \\
& = & \Big\| \Big( \sum_{s \in \sigma'} \mathcal{E}_s|b_s^{\sigma'}|^2 \Big)^{\frac12} \Big\|_s \ \le \ C_s \Big\| \Big( \sum_{s \in\sigma'} |b_s^{\sigma'}|^2 \Big)^{\frac12} \Big\|_{s}.
\end{eqnarray*}
On the other hand, by factorization of conditional expectations \cite[Proposition 2.8]{Jun02} we find an isometric right $\mathcal M_s$-module map $u_s:\mathcal M\rightarrow C(\mathcal M_s)$ for each $s\in\sigma$ and such that $$\mathcal{E}_s(a_s^{\sigma'} b_s^{\sigma'}) = u_s(a_s^{\sigma'*})^*u_s(b_s^{\sigma'}).$$ Given $s \in \sigma'$, let us write $t_s$ in what follows to denote the only member of $\sigma$ which satisfies that $s \in J_{t_s}$. By H\"older inequality and the dual version of Doob's maximal inequality in $L_{w/2}(\mathcal M)$, we find $$\big\| (a_t^\sigma)_{t\in\sigma} \big\|_{L_w(\ell^r_2)} = \Big\| \sum_{t\in\sigma}\sum_{s\in J_t} u_s(a^{\sigma'*}_s)^*u_s(b_s^{\sigma'}) (b_t^\sigma)^{-1} \otimes e_{1,t} \Big\|_w = \big\| A B \big\|_w$$ with
\begin{eqnarray*}
A & = & \sum_{t\in\sigma} \sum_{s\in J_t} u_s(a^{\sigma'*}_s)^* \otimes e_{1,s} \otimes e_{1,t}, \\
B & = & \sum_{t\in\sigma} \sum_{s\in J_t}  u_s(b_s^{\sigma'}) (b_t^\sigma)^{-1} \otimes e_{s,1} \otimes e_{t,t}.
\end{eqnarray*}
This gives rise to the following estimate 
\begin{eqnarray*}
\big\| (a_t^\sigma)_{t\in\sigma} \big\|_{L_w(\ell^r_2)} & \le & \big\| ( A^* A )^\frac12 \big\|_w \ \sup_{t\in\sigma} \Big\| \sum_{s\in J_t}  u_s(b_s^{\sigma'}) (b_t^\sigma)^{-1} \otimes e_{s,1} \Big\|_\infty \\ & = & \Big\| \Big( \sum_{s \in \sigma'} \mathcal{E}_s  |a_s^{\sigma'*}|^2 \Big)^\frac12 \Big\|_w \ \sup_{t\in\sigma} \Big\| (b_t^\sigma)^{-1} \sum_{s\in J_t}  \mathcal{E}_s |b_s^{\sigma'}|^2 (b_t^\sigma)^{-1} \Big\|_\infty^\frac12. 
\end{eqnarray*}
According to the definition of $b_t^\sigma$ and by dual Doob's inequality $$\big\| (a_t^\sigma)_{t\in\sigma} \big\|_{L_w(\ell^r_2)} \le C_w \Big\| \Big( \sum_{s \in \sigma'} |a_s^{\sigma'*}|^2 \Big)^\frac12 \Big\|_w.$$ Altogether, we have proved $\|x\|_{h^{1_c}_{pw}(\sigma)} \le C_{w,s} (\|x\|_{h^{1_c}_{pw}(\sigma')}+\varepsilon)$ for all $\varepsilon > 0$. \fin

\begin{remark} 
\emph{In particular, up to an absolute constant, the above defined norms for column/diagonal Hardy spaces do not depend on the choice of the ultrafilter $\mathcal{U}$. }
\end{remark}

\begin{remark}
\emph{If $D_p$ denotes the optimal constant in the dual form of Doob's inequality over $L_p(\M)$, the proof of Lemma \ref{lem:monotonicity} establishes that the constant $C_{w,s}$ in the statement satisfies $$C_{w,s} \le 2 \sqrt{D_{s/2} D_{w/2}}.$$}
\end{remark}

We are now ready to introduce the algebraic atomic Hardy spaces, which are defined as follows. Let $1\leq p<2$ and $w,s\geq2$ with $1/p=1/w+1/s$. We start with certain auxiliary Hardy spaces given by
\begin{eqnarray*}
\widetilde{h}^c_{pw}(\mathcal M) & = & \big\{ x \in L_p(\mathcal M) : \|x\|_{h^c_{pw}}<\infty \big\}, \\
\widetilde{h}^{1_c}_{pw}(\mathcal M) & = & \big\{ x\in L_p(\mathcal M) : \|x\|_{h^{1_c}_{pw}}<\infty \big\}.
\end{eqnarray*}
The algebraic atomic Hardy spaces are then constructed accordingly 
\begin{eqnarray*}
h^c_{pw}(\mathcal M) & = & \Big\{ \summ_j\lambda_jx_j : x_j\in \widetilde{h}^c_{p_jw_j}\;\mathrm{with}\;w_j>w,\;s_j>s \Big\}, \\
h^{1_c}_{pw}(\mathcal M) & = & \Big\{ \summ_j\lambda_jx_j : x_j\in \widetilde{h}^{1_c}_{p_jw_j}\;\mathrm{with}\;w_j>w,\;s_j>s \Big\},
\end{eqnarray*}
with norm defined via
\begin{eqnarray*}
\|x\|_{h^{c}_{pw}}' & = & \inf \Big\{\summ_j|\lambda_j| : x=\summ_j\lambda_jx_j,\|x_j\|_{h^{c}_{p_jw_j}} \le1 \Big\}, \\
\|x\|_{h^{1_c}_{pw}}' & = & \inf \Big\{\summ_j|\lambda_j| : x=\summ_j\lambda_jx_j,\|x_j\|_{h^{1_c}_{p_jw_j}} \le1 \Big\}.
\end{eqnarray*}
Of course, the row spaces ${h}^{r}_{pw}(\mathcal M)$ and ${h}^{1_r}_{pw}(\mathcal M)$ are defined in a similar way.

\begin{remark}\label{rem:equi norm}
\emph{In the above definitions, summation over $j$ is taken in the $h^c_{pw}$ and $h^{1_c}_{pw}$-norms respectively. Moreover, it is easy to check that the $h^c_{pw}$-norm is equal to its atomic form $$\|x\|_{h_{pw}^c} = \|x\|_{h_{pw}^c}'$$ on the atomic Hardy space $h^c_{pw}$. In particular, we shall drop the $'$ in what follows. All these results hold replacing $c$ by $1_c, r, 1_r$. Moreover, all these spaces are Banach. }
\end{remark}

%\begin{remark}\label{rem:completion}
%Firstly, Secondly, adapting the proof of Lemma 2.27 of \cite{JuPe14}, we can show that $\widetilde{h}^{1_c}_{pw}(\mathcal M)$ is a Banach space due to Lemma %\ref{lem:monotonicity} (ii). Finally, even though we do not know in general whether $\widetilde{h}^{c}_{pw}(\mathcal M)$ is complete, but we can show %$\widetilde{h}^{c}_{p2}%(\mathcal M)$ is complete since 
%$$h^c_{p}(\mathcal M)=\widetilde{h}^{c}_{p2}(\mathcal M).$$
%Indeed, recall that $h^c_{p}(\mathcal M)$ is the completion of $\mathcal M$ under the norm $\|\cdot\|_{h^c_p}:=\lim_{\sigma,\mathcal U}\|\cdot\|_{h^c_p(\sigma)}$, hence it is %trivial that $h^c_{p}(\mathcal M)\subset \widetilde{h}^{c}_{p2}(\mathcal M);$ On the other hand, for any $x\in \widetilde{h}^{c}_{p}$, we can find a partition such that
%$$\|x\|_{{h}^{c}_{p}(\sigma)}\leq 2\|x\|_{{h}^{c}_{p}}<\infty.$$ By the fact $\M$ is dense in ${{h}^{c}_{p}(\sigma)}$, for any $\varepsilon>0$, we find $y\in\M$ such that $\|x-y\|%_{{h}^{c}_{p}(\sigma)}<\varepsilon$. Therefore by the decreasing property of the ${{h}^{c}_{p}(\sigma)}$-norms---Lemma 4.13 of \cite{JuPe14},
%$$\|x-y\|_{{h}^{c}_{p}}\leq \|x-y\|_{{h}^{c}_{p}(\sigma)}<\varepsilon,$$ which implies $x\in h^{c}_{p}(\mathcal M)$.
%\end{remark}

The proof of Theorem Biii crucially rests on the algebraic atomic characterization in Theorem \ref{pro:hat hcpw=hcpw} below. Its proof is quite involved, so that we postpone it for the next subsection and complete the proof of Theorem Biii taking it for granted.  

\begin{theorem}\label{pro:hat hcpw=hcpw}
We have $$h^c_{pw}(\mathcal M) \simeq \widehat{h}^c_{pw}(\mathcal M) \quad \mbox{and} \quad h^{1_c}_{pw}(\mathcal M) = \widehat{h}^{1_c}_{pw}(\mathcal M)$$
for $1\leq p<2$ and $w,s\geq2$ with $\frac1p = \frac1w + \frac1s$. The same holds true for row spaces.
\end{theorem}

\noindent \textbf{Proof of Theorem Biii.}
Let $1<p<2$ and $w>2$. By Theorem \ref{pro:hat hcpw=hcpw},
it suffices to show that $h^c_{pw}(\mathcal M) \subset \mathcal H^c_p(\mathcal M)$ and $h^{1_c}_{pw}(\mathcal M)\subset\mathcal H^c_p(\mathcal M)$ respectively. For the  first inclusion, given $x \in h^c_{pw}(\mathcal M)$ we may assume that $$x \in \widetilde{h}^c_{\widetilde{p}\widetilde{w}}(\mathcal M) \quad \mbox{for some} \quad \widetilde{p}>p, \widetilde{w}>w$$ with $h_{\widetilde{p} \widetilde{w}}^c$-norm smaller than 1. By Lemma \ref{lem:monotonicity}, there exists one partition $\sigma$ such that $\|x\|_{{h}^c_{\widetilde{p}\widetilde{w}}(\sigma)}<1$. Using Davis decomposition \cite[Theorem 6.3.2]{JuPe14}, the decreasing property of ${{h}^{c}_{p}(\sigma)}$-norms  \cite[Lemma 5.3.1]{JuPe14} and arguing as in \cite[Lemma 2.9]{HJP}, we conclude $$\|x\|_{\mathcal H^c_p}\leq \|x\|_{h^c_p}\leq\|x\|_{h^c_p(\sigma)} \le c_{\widetilde{w}}\|x\|_{{h}^c_{\widetilde{p}\widetilde{w}}(\sigma)} < c_{\widetilde{w}}.$$ The second inclusion is proved similarly. Take $x\in \widetilde{h}^{1_c}_{\widetilde{p}\widetilde{w}}(\mathcal M)$ for some $\widetilde{p}>p, \widetilde{w}>w$ with norm smaller than 1. Given any partition $\sigma$ and arguing once more as in the proof of \cite[Lemma 2.9]{HJP} yields |in conjunction with Lemma \ref{lem:monotonicity} above| the following estimate $$\|x\|_{H^c_p(\sigma)}\lesssim \|x\|_{h^{1_c}_{\widetilde{p}\widetilde{w}}(\sigma)} < 1.$$
Hence, taking suprema over all $\sigma \in \Sigma$ $$\|x\|_{\mathcal{H}^c_p} \sim \sup_{\sigma \in \Sigma} \|x\|_{H^c_p(\sigma)} \lesssim 1$$ where we use the increasing property of the ${{H}^{c}_{p}(\sigma)}$-norms \cite[Lemma 3.3.1]{JuPe14}. \fin

\subsection{Atomic decomposition.}

We finsih the paper with the proof of Theorem \ref{pro:hat hcpw=hcpw}: 

\noindent \textbf{The inclusions $h^c_{pw}(\mathcal M)\subset\widehat{h}^c_{pw}(\mathcal M)$ and $h^{1_c}_{pw}(\mathcal M)\subset\widehat{h}^{1_c}_{pw}(\mathcal M)$.} We only deal with $h^c_{pw}(\mathcal M)$, since similar arguments work also for the diagonal Hardy spaces. By definition, given $x\in h^c_{pw}(\mathcal M)$ we may find $p<p_j<2$, $w_j>w$ and $s_j>s$ satisfying $x = \sum_{j\geq1}\lambda_jx_j$ for some $x_j$ with $\|x_j\|_{h^c_{p_jw_j}} \le 1$ and $$\sum_{j\geq1}|\lambda_j|\leq (1+\delta) \|x\|_{h^c_{pw}}.$$ On the other hand, given any $\sigma \in \Sigma$ we may write in turn $$x_j=\sum_{t\in\sigma}a^\sigma_{t}(j)b^\sigma_{t}(j)$$ with $\mathcal{E}_{t}(a^\sigma_{t}(j))=0$, $b^\sigma_{t}(j)\in L_{s_j}(\mathcal{M}_{t})$ for all $j \ge 1$ and such that $$\big\| (a^\sigma_{t}(j))_{t\in\sigma} \big\|_{L_{w_j}(\ell^r_2)} = \big\| (b^\sigma_{t}(j))_{t\in\sigma} \big\|_{L_{s_j}(\ell^c_2)} \le (1 + \delta) \|x_j\|_{h_{p_jw_j}^c(\sigma)}^\frac12.$$ Define $$a_t^\sigma = \sum_{j \ge 1} {\lambda_j}a^\sigma_t(j)b^\sigma_t(j)(b^\sigma_t)^{-1} \quad \mbox{where} \quad b_t^\sigma = \Big( \sum_{j \ge 1} |\lambda_j| |b^\sigma_t(j)|^2 \Big)^{\frac12},$$ where we assuming $b_t^\sigma$ invertible by approximation. This gives $$x = \sum_{j\geq1}\lambda_jx_j = \sum_{t\in\sigma}a^\sigma_{t}b^\sigma_{t}$$ for all $\sigma \in \Sigma$. Therefore, letting $a = (a^\sigma)^\bullet$ and $b = (b^\sigma)^\bullet$ we obtain $$\|x\|_{\widehat{h}_{pw}^c} \le \|a^*\|_{X_w^\circ} \|b\|_{X_s^{ad}} = \lim_{\sigma, \mathcal{U}} \Big\| \Big( \sum_{t \in \sigma} a_t^\sigma a_t^{\sigma*} \Big)^\frac12 \Big\|_w \Big\| \Big( \sum_{t \in \sigma} b_t^{\sigma*} b_t^{\sigma} \Big)^\frac12 \Big\|_s \ = \ \lim_{\sigma,\mathcal{U}} A_\sigma B_\sigma.$$ Indeed, considering finite $j$-truncations of $a_t^\sigma$ and $b_t^\sigma$ it is easy to check that $a^* \in X_w^\circ$ and $b \in X_s^{ad}$. Once we have estimated the norm of $x$ by the ultralimit of $A_\sigma B_\sigma$, it suffices to show the following estimates 
\begin{eqnarray*}
A_\sigma & \le & \Big( \sum_{j \ge 1} |\lambda_j| \big\| (a^\sigma_{t}(j))_{t\in\sigma} \big\|_{L_{w_j}(\ell^r_2)}^2 \Big)^\frac12, \\
B_\sigma & \le & \Big( \sum_{j \ge 1} |\lambda_j| \hskip1pt \big\| (b^\sigma_{t}(j))_{t\in\sigma} \big\|_{L_{s_j}(\ell^c_2)}^2 \hskip1pt \Big)^\frac12, 
\end{eqnarray*}
since then, taking ultralimits, we obtain the expected estimate $$\|x\|_{\widehat{h}_{pw}^c} \le (1+\delta)^2 \sum_{j \ge 1} |\lambda_j| \le (1+\delta)^3 \|x\|_{h^c_{pw}}$$ since $\delta > 0$ may be taken arbitrarily small. To justify the claim for $A_\sigma$ $$a_t^\sigma = \Big( \sum_{j \ge 1} |\lambda_j|^\frac12 a^\sigma_t(j) \otimes e_{1j} \Big) \Big( \sum_{j \ge 1} \lambda_j |\lambda_j|^{-\frac12} b^\sigma_t(j)(b^\sigma_t)^{-1} \otimes e_{j1} \Big) = \alpha_t^\sigma \beta_t^\sigma,$$ where $\beta_t^\sigma$ is clearly a contraction, so that $a_t^\sigma a_t^{\sigma*} \le \alpha_t^\sigma \alpha_t^{\sigma*}$. This implies that
\begin{eqnarray*}
A_\sigma & \le & \Big\| \sum_{t \in \sigma} \alpha_t^\sigma \alpha_t^{\sigma*} \Big\|_{\frac{w}{2}}^\frac12 \\ & = & \Big\| \sum_{j \ge 1} |\lambda_j| \sum_{t \in \sigma} a_t^\sigma(j) a_t^{\sigma}(j)^* \Big\|_{\frac{w}{2}}^\frac12 \\ & \le & \Big( \sum_{j \ge 1} |\lambda_j| \Big\| \sum_{t \in \sigma} a_t^\sigma(j) a_t^{\sigma}(j)^* \Big\|_{\frac{w}{2}} \Big)^\frac12
\end{eqnarray*}
since $w_j \ge w$. The claim for $B_\sigma$ follows by the triangle inequality in $L_{s_j/2}(\M)$. \fin 

\noindent \textbf{The inclusion $\widehat{h}^{1_c}_{pw}(\mathcal M)\subset h^{1_c}_{pw}(\mathcal M)$.}  Given $x\in \widehat{h}^{1_c}_{pw}(\mathcal M)$ and $\delta>0$, there exists a representation $x=\langle\langle a^*,b\rangle\rangle$ such that $\|a^*\|_{X_w}\|b\|_{X_s}\leq (1+\delta)\|x\|_{\widehat{h}^{1_c}_{pw}}$. By density we can further assume that $$a\in \prod_{\sigma,\mathcal{U}}L_{\widetilde{w}}\big(\M;\ell^r_2(\sigma\times\mathbb N)\big) \quad \mbox{and} \quad b\in \prod_{\sigma,\mathcal{U}}L_{\widetilde{s}}\big(\M;\ell^c_2(\sigma\times\mathbb N)\big)$$ for some $\widetilde{w}>w$ and $\widetilde{s}>s$ such that $$\|a^*\|_{X_{\widetilde{w}}}\leq (1+\delta)\|a^*\|_{X_w} \hskip16pt \mbox{and} \hskip24pt \|b\|_{X_{\widetilde{s}}}\leq (1+\delta)\|b\|_{X_s}.$$ We claim that it suffices to prove the following estimate 
\begin{equation} \label{w limit controlled}
\|x\|_{h^{1_c}_{\widetilde{p}\widetilde{w}}}\leq \lim_{\sigma,\mathcal{U}} \Big\| \underbrace{\sum_{t\in\sigma}d^\sigma_{t}(a^\sigma_{t}b^\sigma_{t})}_{x_\sigma} \Big\|_{h^{1_c}_{\widetilde{p}\widetilde{w}}(\sigma)}.
\end{equation}
Indeed, in that case we conclude $$\|x\|_{h_{pw}^{1_c}} = \lim_{\sigma,\mathcal{U}} \|x\|_{h_{pw}^{1_c}(\sigma)} \le \lim_{\sigma,\mathcal{U}} \|x\|_{h_{\widetilde{p}\widetilde{w}}^{1_c}(\sigma)} \le \lim_{\sigma,\mathcal{U}} \|x_\sigma\|_{h_{\widetilde{p}\widetilde{w}}^{1_c}(\sigma)} \le (1+\delta)^3 \|x\|_{\widehat{h}_{pw}^{1_c}}$$ with $\delta > 0$ arbitrarily small. Let us then prove the claim. According to the definition of $\widehat{h}_{\widetilde{p}\widetilde{w}}^{1_c}(\M)$, we know that $x = w-L_{\widetilde{p}}-\lim_{\sigma,\mathcal{U}} x_\sigma$. By the monotonicity Lemma \ref{lem:monotonicity}, it suffices to prove $$\|x\|_{h^{1_c}_{\widetilde{p}\widetilde{w}}(\sigma_0)}\leq\lim_{\sigma,\mathcal{U}}\|x_\sigma\|_{{h^{1_c}_{\widetilde{p}\widetilde{w}}}(\sigma)} \quad \mbox{for all $\sigma_0 \in \Sigma$}.$$ Let us fix a $\sigma_0 \in \Sigma$ and $\varepsilon>0$. According to Mazur's lemma, we can find a finite sequence of positive numbers $(\alpha_m)^M_{m=1}$ with $\sum_m \alpha_m=1$ and partitions $\sigma^1, \sigma^2, \ldots, \sigma^M$ containing $\sigma_0$ such that 
\begin{equation} \label{mazur}
\Big\| x - \sum_{m=1}^M \alpha_m x_{\sigma^m} \Big\|_{\widetilde{p}} < \varepsilon \quad \mbox{and} \quad \|x_{\sigma^m}\|_{h^{1_c}_{\widetilde{p}\widetilde{w}}(\sigma^m)} \le (1+\varepsilon) \lim_{\sigma,\mathcal{U}} \|x_\sigma\|_{h^{1_c}_{\widetilde{p}\widetilde{w}}(\sigma)}
\end{equation}
for all $m=1,2, \ldots, M$. Hence, we may decompose $$\|x\|_{h^{1_c}_{\widetilde{p}\widetilde{w}}(\sigma_0)} \le \Big\| x - \sum_{m=1}^M\alpha_mx_{\sigma^m} \Big\|_{h^{1_c}_{\widetilde{p}\widetilde{w}}(\sigma_0)} + \Big\| \sum_{m=1}^M \alpha_m x_{\sigma^m} \Big\|_{h^{1_c}_{\widetilde{p}\widetilde{w}}(\sigma_0)}.$$ The second quantity on the right hand side is smaller than $$(1+\varepsilon)\lim_{\sigma,\mathcal{U}}\|x_\sigma\|_{h^{1_c}_{\widetilde{p}\widetilde{w}}(\sigma)}$$ by the triangle inequality and Lemma \ref{lem:monotonicity} |since $\sigma_{m}$ contains $\sigma_0$| together with the second estimate in \eqref{mazur}. On the other hand, the first quantity can be made arbitrarily small as a consequence of the first estimate in \eqref{mazur} and the inequality below $$\|z\|_{h^{1_c}_{\widetilde{p}\widetilde{w}}(\sigma_0)}\leq |\sigma_0|\|z\|_{\widetilde{p}},$$ which is valid for $1<{\widetilde{p}}<2$ and any finite partition $\sigma_0$. Indeed 
\begin{eqnarray*}
\hskip40pt \|z\|_{h^{1_c}_{\widetilde{p}\widetilde{w}}(\sigma_0)} & = & \Big\| \sum_{t\in\sigma_0} d^{\sigma_0}_t (z|z|^{-\widetilde{p}/\widetilde{s}}|z|^{\widetilde{p}/\widetilde{s}}) \Big\|_{h^{1_c}_{\widetilde{p}\widetilde{w}}(\sigma_0)} \\ & \le & |\sigma_0|^{\frac12} \big\| z|z|^{-\widetilde{p}/\widetilde{s}} \big\|_{\widetilde{w}} |\sigma_0|^\frac12 \big\| |z|^{\widetilde{p}/\widetilde{s}} \big\|_{\widetilde{s}} \ \le \ |\sigma_0| \|z\|_{\widetilde{p}}. \hskip45pt \square
\end{eqnarray*}

The inclusion $\widehat{h}^c_{pw}(\mathcal M)\subset h^c_{pw}(\mathcal M)$ requires different arguments than those which we have used for diagonal Hardy spaces, since the column Hardy space norms have opposite monotonicity and we can not prove something similar to \eqref{w limit controlled}. We shall proceed with the proof after a couple of preliminary lemmas.

\begin{lemma}\label{lem:two approximation}
Let $x\in \widehat{h}^c_{pw}(\mathcal M)$ and $\varepsilon, \varepsilon' > 0$, then$\hskip1pt :$ 
\begin{itemize}
\item[i)] Given $\widetilde{w}<\infty$, there exists $x'\in \widehat{h}^c_{pw}(\mathcal M) \cap L_{\widetilde{w}}(\mathcal M)$ with $\|x-x'\|_{\widehat{h}^c_{pw}}<\varepsilon$.

\item[ii)] There exists $x''\in \widetilde{h}^c_{pw}(\mathcal M)$ with $\|x'-x''\|_p\leq\varepsilon'$ and $\|x''\|_{h^c_{pw}}\leq 4\|x'\|_{\widehat{h}^c_{pw}}$.
\end{itemize}
\end{lemma}

\dem Since $x \in \widehat{h}^c_{pw}(\mathcal M)$, there exists a representation $$x=\langle a^*,b\rangle \quad \mbox{such that} \quad \max \big\{ \|a^*\|_{X^\circ_w}, \|b\|_{X^{ad}_s} \big\} \le 2 \|x\|^{\frac12}_{\widehat{h}^c_{pw}}.$$ Let $v>\max(2\widetilde{w},w,s)$ and $\delta > 0$. By the density we may find $$\widetilde{a}^*, \widetilde{b}\in \prod_{\sigma,\mathcal{U}}L_v(\mathcal M;\ell^c_2(\sigma))$$ such that $$\max \big\{ \| a^*-\widetilde{a}^* \|_{X^\circ_{w}}, \| b - \widetilde{b} \|_{X^{ad}_s} \big\} \le \delta.$$ Then we define $x'=\langle \hskip1pt \widetilde{a}^*,\widetilde{b} \rangle$ and it clearly has the desired properties 
\begin{itemize}
\item[$\bullet$] $x'\in \widehat{h}^c_{pw}(\mathcal M)\cap L_{\widetilde{w}}(\mathcal M)$, 

\item[$\bullet$] If we pick $\delta = \sqrt{4  \|x\|_{\widehat{h}^c_{pw}} + \varepsilon} - 2  \|x\|^{\frac{1}{2}}_{\widehat{h}^c_{pw}}$, then 
\begin{eqnarray*}
\|x-x'\|_{\widehat{h}^c_{pw}} &\le & \|a^*-\widetilde{a}^*\|_{X^\circ_w}\|b\|_{X^{ad}_s} \\ & + & \|\widetilde{a}^* -a^*\|_{X^\circ_w} \|b-\widetilde{b}\|_{X^{ad}_s} \\ & + & \|a^*\|_{X^\circ_w} \|b-\widetilde{b}\|_{X^{ad}_s} \ \le \ 4\delta \|x\|^{\frac{1}{2}}_{\widehat{h}^c_{pw}} + \delta^2 \ = \ \varepsilon.
\end{eqnarray*}
\end{itemize}
Let us now prove the second assertion. Once we know that $x' \in \widehat{h}_{pw}^c(\M)$ we know that $$x' = \big\langle {\widehat{a}}^*,{\widehat{b}} \big\rangle = w - L_p - \lim_{\sigma,\mathcal{U}} \sum_{t \in \sigma} {\widehat{a}^\sigma}_t {\widehat{b}^\sigma}_t \quad \mbox{with} \quad \|{\widehat{a}}^* \|_{X_w^\circ} \|\widehat{b}\|_{X_s^{ad}}\leq 2\|x'\|_{\widehat{h}^c_{pw}}.$$ Define $x'_\sigma=\sum_{t\in\sigma}\widehat{a}^\sigma_{t}\widehat{b}^\sigma_{t}$. As we did in \eqref{mazur} we can find a sequence of positive numbers $(\alpha_m)^M_{m=1}$ satisfying $\sum_{m}\alpha_m=1$ and also partitions $\sigma^1, \sigma^2, \ldots, \sigma^M$ such that
$$\Big\| x' - \sum_{m=1}^M \alpha_m x'_{\sigma^m} \Big\|_p < \varepsilon' \quad \mbox{and} \quad \|x'_{\sigma^m}\|_{h^c_{pw}(\sigma^m)} \le 2 \lim_{\sigma,\mathcal{U}} \|x'_\sigma\|_{h^{c}_{pw}(\sigma)} \le 4\|x'\|_{\widehat{h}^c_{pw}}$$ for all $m=1,2, \ldots, M$. Defining $x''=\sum_m\alpha_mx_{\sigma^m}$ and applying Lemma \ref{lem:monotonicity}, we conclude \\ [3pt] \null \hskip35pt $ \displaystyle \|x''\|_{h^{c}_{pw}}\leq\sum_{m=1}^M\alpha_m\|x'_{\sigma^m}\|_{h^c_{pw}}\leq\sum_{m=1}^M\alpha_m\|x'_{\sigma^m}\|_{h^c_{pw}(\sigma^m)}\leq 4\|x'\|_{\widehat{h}^c_{pw}}.$ \hfill \fin

\begin{lemma}\label{lem:equiv norm on subspace}
Let $x\in h^c_{pw}(\mathcal M)$ such that for any $\varepsilon>0$, there exist $x(\varepsilon)\in \widetilde{h}^c_{pw}(\mathcal M)$ satisfying the inequalities $\|x-x(\varepsilon)\|_p\leq\varepsilon$ and $\|x(\varepsilon)\|_{h^c_{pw}}\leq K\|x\|_{\widehat{h}^c_{pw}}$. Then we find $$\|x\|_{h^c_{pw}}\leq K\|x\|_{\widehat{h}^c_{pw}}.$$
\end{lemma}

\dem Since $h^c_{pw}(\mathcal M)$ is injective in $L_p(\mathcal M)$, it turns out that $L_{p'}(\mathcal M) = L_p(\M)^*$ is dense in the dual space $h^c_{pw}(\mathcal M)^*$. Given any $z \in L_{p'}(\mathcal M)$, we have the following estimate
\begin{eqnarray*}
|\varphi_x(z)| & = & |\tau(x^*z)| \ = \ \big| \lim_{\varepsilon\rightarrow0}\tau(x(\varepsilon)^*y) \big| \\
&\leq & \lim_{\varepsilon\rightarrow0} \|x(\varepsilon)\|_{h^c_{pw}} \|z\|_{(h^c_{pw})^*} \ \le \ K\|x\|_{\widehat{h}^c_{pw}} \|z\|_{(h^c_{pw})^*}
\end{eqnarray*}
which implies the desired estimate by taking supreme over $z \in L_{p'}(\M)$. \fin

\noindent \textbf{The inclusion $\widehat{h}^c_{pw}(\mathcal M)\subset h^c_{pw}(\mathcal M)$.} Let $x \in \widehat{h}^c_{pw}(\mathcal M)$. By the first part of Lemma \ref{lem:two approximation}, for any $\varepsilon>0$ and $\widetilde{w} < \infty$, we find $x(\varepsilon)\in \widehat{h}^c_{pw}(\mathcal M)\cap L_{\widetilde{w}}(\mathcal M)$. Actually $x(\varepsilon)$ can be chosen to be in $\widetilde{h}^{c}_{\widetilde{p}\widetilde{w}}(\mathcal M)$ where $p<\widetilde{p}<2$, $\widetilde{w}>w$, $\widetilde{s}>2$ with $$\frac{1}{\widetilde{p}} = \frac{1}{\widetilde{w}} + \frac{1}{\widetilde{s}}.$$ Indeed, fix $\sigma \in \Sigma$ and write $x(\varepsilon)=\sum_{t\in\sigma}d^\sigma_{t}(x(\varepsilon))\cdot \1$. By monotonicity Lemma \ref{lem:monotonicity} $$\|x(\varepsilon)\|_{h^c_{\widetilde{p}\widetilde{w}}} \le \|x(\varepsilon)\|_{h^c_{\widetilde{p}\widetilde{w}}(\sigma)} \le \Big\|\Big( \sum_{t\in\sigma} |d^\sigma_{t}(x(\varepsilon)^*)|^2 \Big)^{\frac{1}{2}} \Big\|_{\widetilde{w}} |\sigma|^{\frac{1}{2}} \le C_{\widetilde{w}} |\sigma|^{\frac12} \|x(\varepsilon)\|_{\widetilde{w}}.$$ The last estimate follows from the noncommutative form of Burkholder-Gundy inequality \cite{PX}. Now, we know from the second part of Lemma \ref{lem:two approximation} that $x(\varepsilon)$ satisfies the conditions in Lemma \ref{lem:equiv norm on subspace} with $K=4$ for any $\varepsilon>0$. Next, if we apply Lemma \ref{lem:two approximation} for all $\varepsilon_k =2^{-k}\|x\|_{\widehat{h}^c_{pw}}$, we obtain $$x=x(\varepsilon_1)+\sum_{k\geq1}(x(\varepsilon_{k+1})-x(\varepsilon_k))$$ where the summation converges in $\widehat{h}^c_{pw}(\mathcal M)$. Recall that all $x(\varepsilon_{k+1})-x(\varepsilon_k)$'s satisfy the conditions in Lemma \ref{lem:equiv norm on subspace} with $K=8$. In particular, putting altogether the assertion follows 
\begin{eqnarray*}
\|x\|_{{h}^c_{pw}} & \le & \|x(\varepsilon_1)\|_{h^c_{pw}} + \sum_{k \ge 1} \| x(\varepsilon_{k+1})-x(\varepsilon_k) \|_{h^c_{pw}} \\ & \lesssim & \|x(\varepsilon_1)\|_{\widehat{h}^c_{pw}} + \sum_{k \ge 1} \| x(\varepsilon_{k+1})-x(\varepsilon_k) \|_{\widehat{h}^c_{pw}} \ \lesssim \ \Big( 1 + \sum_{k \ge 1} 2^{-k} \Big) \|x\|_{\widehat{h}^c_{pw}}.
\end{eqnarray*}
This completes the proof of $\widehat{h}^c_{pw}(\mathcal M)\subset h^c_{pw}(\mathcal M)$ and thus of Theorem \ref{pro:hat hcpw=hcpw}. \fin

\noindent \textbf{Acknowledgement.} Junge is partially supported by the NSF DMS-1201886 and NSF DMS-1501103. Parcet is partially supported by Proyecto Intramural 201650E030 (CSIC) and Proyecto Excelencia Europa QHA MTM2016-81700-ERC (MINECO). Hong is partially supported by the NSF of China-11601396, Funds for Talents of China-413100002 and 1000 Young Talent Researcher Programm of China-429900018-101150(2016).

\vskip-3pt

\hfill \noindent \textbf{Guixiang Hong} \\
\null \hfill School of Mathematics and Statistics \\ 
\null \hfill Wuhan University \\ 
\null \hfill Wuhan 430072. China \\
\null \hfill\texttt{guixiang.hong@whu.edu.cn}

\

\hfill \noindent \textbf{Marius Junge}\\
\null \hfill Department of Mathematics\\ \null \hfill
University of Illinois at Urbana-Champaign\\ \null \hfill
1409 W. Green St. Urbana, IL 61891. USA \\ \null \hfill
\texttt{junge@math.uiuc.edu}

\

\hfill \noindent \textbf{Javier Parcet} \\
\null \hfill Instituto de Ciencias Matem{\'a}ticas \\ \null \hfill
CSIC-UAM-UC3M-UCM \\ \null \hfill Consejo Superior de
Investigaciones Cient{\'\i}ficas \\ \null \hfill C/ Nicol\'as Cabrera 13-15.
28049, Madrid. Spain \\ \null \hfill\texttt{javier.parcet@icmat.es}

\end{document}